\documentstyle[12pt,bezier,epsfig]{article}

\newcommand{\Exp}{{\rm I\hspace{-0.8mm}E}}
\newcommand{\Prob}{{\rm I\hspace{-0.8mm}P}}
\newcommand{\Var}{{\bf Var\,}}
\newcommand{\Cov}{{\bf Cov}}

\newcommand{\iz}{{\rm \rlap Z\kern 2.2pt Z}}

\newcommand{\RL}{{\rm I\hspace{-0.8mm}R}}

\newcommand{\bbN}{{\rm I\hspace{-0.8mm}N}}

\newcommand{\proof}{\noindent {\bf Proof:} \ }
\newcommand{\halmos}{\newline\vspace{3mm}\hfill $\Box$}

\newtheorem{theorem}{Theorem}
\newtheorem{lemma}{Lemma}

\newtheorem{corollary}{Corollary}

\newtheorem{remark}{Remark}

\title{\bf Remarks on Pickands theorem}
\author{
\\
ZBIGNIEW MICHNA\\
Department of Mathematics and Cybernetics\\
Wroc{\l}aw University of Economics\\
Wroc{\l}aw, Poland}
\date{}

\begin{document}

\maketitle

\begin{abstract}
In this article we present Pickands theorem and his double sum method. We follow Piterbarg's
proof of this theorem. 
Since his proof relies on general lemmas we present a complete proof of Pickands
theorem using Borell inequality and Slepian lemma.
The original Pickands proof is rather complicated
and is mixed with upcrossing probabilities for stationary Gaussian processes. We give a lower bound for 
Pickands constant.

\vspace{5mm}
\newline
\vspace{2cm}
\end{abstract}

\section{Introduction}
James Pickands III (see \cite{pi:69a} and \cite{pi:69b}) gave an elegant and sophisticated way of finding the asymptotic behavior of the probability 
$$ 
\Prob(\sup_{t\in {\bf T}}X(t)>u)
$$
as $u\rightarrow\infty$ where $X$ is a Gaussian process. More precisely for $t\in [0, p]$ let $X(t)$ be a continuous stationary Gaussian 
process with expected value $\Exp X(t)=0$ and covariance 
$$
r(t)=\Exp(X(t+s)X(s))=1-|t|^\alpha+o(|t|^\alpha)
$$
where $0<\alpha\leq 2$.
Furthermore we assume that $r(t)<1$ for all $t>0$. Then
$$
\Prob(\sup_{t\in [0,p]} X(t)>u)=H_\alpha\, p\, u^{2/\alpha}\,\Psi(u)(1+o(1))
$$
where $H_{\alpha}$ is a positive and finite constant (Pickands constant) and $\Psi(u)$ is the tail of the standard normal distribution.
We will follow Piterbarg's
proof of this theorem. Since his proof relies on general lemmas we present a complete proof of Pickands
theorem using Borel inequality and Slepian lemma. 
Lemma \ref{lemdlak} below is different than Lemma D.2. in Piterbarg \cite{pit:96} that is the constant before 
exponent depends on $T$.

The original Pickands proof is rather complicated
and is mixed with upcrossing probabilities for Gaussian stationary processes. In his paper this
theorem is a lemma (see \cite{pi:69b}). The proof of Pickands theorem is based on the elementary Bonferroni
inequality which in the literature is in a too strong version. In this paper we present a sharper version
of the Bonferroni inequality which has an impact on some lower bounds of Pickands constant (see \cite{de:mi:ro:03} and \cite{sh:96}). Some upper estimates of Pickands constant can be found in \cite{de:ki:08}.

\section{Lemmas and auxiliary theorems}
In the paper we will consider real-valued stochastic processes and fields.
Let us denote
$$
\Psi(u)=1-\Phi(u)=\frac{1}{\sqrt{2\pi}}\int_u^{\infty}e^{-\frac{s^2}{2}}\,ds
$$
and notice
\begin{equation}\label{glasym}
\Psi(u)=\frac{1}{\sqrt{2\pi}u}\,e^{-\frac{u^2}{2}}\,(1+o(1))
\end{equation}
as $u\rightarrow\infty$. More precisely for $u>0$
$$
\left(\frac{1}{u}-\frac{1}{u^3}\right)\frac{1}{\sqrt{2\pi}}
\,e^{-\frac{u^2}{2}}< \Psi(u)<\frac{1}{u}\,\frac{1}{\sqrt{2\pi}}
\,e^{-\frac{u^2}{2}}\,.
$$
\begin{lemma}\label{gausv}
Let $(X_1,X_2)$ be a Gaussian vector with values in $\RL^2$ with $\Exp X_1=m_1$, $\Exp X_2=m_2$,
$\Var X_1=\sigma_1^2$, $\Var X_2=\sigma_2^2$ and $\rho=\Cov(X_1,X_2)$.
Then
$$
X_2=\alpha X_1+Z
$$
where $$
\alpha=\frac{\rho}{\sigma_1^2}
$$
and $Z$ is independent of $X_1$ and is normally distributed with mean $m_2-\alpha m_1$ and
variance 
$$
\sigma_2^2-\frac{\rho^2}{\sigma_1^2}\,.
$$
\end{lemma}

\begin{lemma}\label{bonf}
{\em (Bonferroni inequality)}
Let $(\Omega, {\cal S}, \Prob)$ be a probability
space and \\
$A_1,A_2,\ldots,A_n\in {\cal S}$ for $n\geq 2$. Then
$$
\Prob(\bigcup_{i=1}^n A_i)\geq \sum_{i=1}^n\Prob(A_i)-
\sum_{1\leq i<j\leq n}\Prob(A_i\cap A_j)\,.
$$
\end{lemma}
\proof
Our proof will follow by induction.
For $n=2$ we have $\Prob(A_1\cup A_2)=\Prob(A_1)+\Prob(A_2)-\Prob(A_1\cap A_2)$.
Thus let us assume that the inequality is true for $n$. Then
\begin{eqnarray*}
\Prob(\bigcup_{i=1}^{n+1} A_i)&=& \Prob(\bigcup_{i=1}^n A_i)+\Prob(A_{n+1})
-\Prob((\bigcup_{i=1}^n A_i)\cap A_{n+1})\\
&=& \Prob(\bigcup_{i=1}^n A_i)+\Prob(A_{n+1})
-\Prob(\bigcup_{i=1}^n (A_i\cap A_{n+1}))\\
&\geq & \sum_{i=1}^{n+1}\Prob(A_i)-\sum_{1\leq i<j\leq n}\Prob(A_i\cap A_j)
-\Prob(\bigcup_{i=1}^n (A_i\cap A_{n+1}))\\
&\geq&\sum_{i=1}^{n+1}\Prob(A_i)-\sum_{1\leq i<j\leq n}\Prob(A_i\cap A_j)
-\sum_{i=1}^n\Prob(A_i\cap A_{n+1})\\
&=&\sum_{i=1}^{n+1}\Prob(A_i)-\sum_{1\leq i<j\leq n+1}\Prob(A_i\cap A_j)
\end{eqnarray*}
where in the third line we used the induction hypothesis. Thus by induction the inequality is
valid for all $n\geq 2$.
\halmos

Using above Bonferroni inequality we get a sharper lower bound of Pickands constant than in \cite{de:mi:ro:03} (twice as big) whose the proof goes the same way as in \cite{de:mi:ro:03}. 
\begin{theorem}
$$
H_\alpha\geq \frac{\alpha}{2^{2+\frac{2}{\alpha}}\Gamma\left(\frac{1}{\alpha}\right)}\,.
$$
\end{theorem}

The next theorem is also elementary but very useful.
\begin{theorem}
{\em (Slepian inequality)}
Let Gaussian fields $X(t)$ and $Y(t)$ be separable where $t\in {\bf T}$ and ${\bf T}$ is
an arbitrary parameter set. Moreover we assume that the covariance functions
$r_X(t,s)=\Exp(X(t)-\Exp X(t))(X(s)-\Exp X(s))$ and \\
$r_Y(t,s)=\Exp(Y(t)-\Exp Y(t))(Y(s)-\Exp Y(s))$
satisfy
$$
r_X(t,t)=r_Y(t,t)\,
$$
$$
r_X(t,s)\leq r_Y(t,s)
$$
for all $t,s\in {\bf T}$ and their expected values fulfill
$$
\Exp X(t)=\Exp Y(t)
$$
for all $t\in {\bf T}$. Then for any $u$
$$
\Prob(\sup_{t\in {\bf T}} X_t<u)\leq 
\Prob(\sup_{t\in {\bf T}} Y_t<u)\,.
$$
\end{theorem}

The next theorem is the most important tool in the theory of Gaussian processes (see \cite{ad:ta:04}).
\begin{theorem}\label{borel}
{\em (Borell inequality)}
Let $X(t)$ be a centered a.s. bounded Gaussian field where $t\in {\bf T}$ and ${\bf T}$ is
an arbitrary parameter set. Then
$$
\Exp\sup_{t\in {\bf T}}X(t)=m< \infty\,,\,\,\,\,\,\,\,\,
\sup_{t\in {\bf T}}\Var X(t)=\sigma^2< \infty\,,
$$
and for all $w\geq m$
$$
\Prob(\sup_{t\in {\bf T}}X(t)>w)\leq \exp\left(-\frac{(w-m)^2}{2\sigma^2}\right)\,.
$$
\end{theorem}
We will assume that $0<\alpha\leq 2$.
The next lemma one can find in Piterbarg \cite{pit:96} but it is in a more general setting which is not necessary
in the proof of Pickands theorem.
\begin{lemma}\label{lemobok}
Let $\chi(t)$ be a continuous Gaussian field where $t=(t_1,t_2)\in\RL^2$ with 
$\Exp\chi(t)=-|t_1|^\alpha-|t_2|^\alpha$ and $\Cov(\chi(t),\chi(s))=|t_1|^\alpha+|t_2|^\alpha+
|s_1|^\alpha+|s_2|^\alpha-|t_1-s_1|^\alpha-|t_2-s_2|^\alpha$ ($s=(s_1,s_2)$) and $X(t)$ be a continuous homogeneous Gaussian 
field where $t=(t_1,t_2)\in\RL^2$ with expected value $\Exp X(t)=0$ and covariance 
$$
r(t)=\Exp(X(t+s)X(s))=1-|t_1|^\alpha-|t_2|^\alpha+o(|t_1|^\alpha+|t_2|^\alpha)\,.
$$
Then for any compact set ${\bf T}\subset\RL^2$
$$
\Prob(\sup_{t\in\, u^{-2/\alpha}{\bf T}}X(t)>u)=\Psi(u)H({\bf T})(1+o(1))
$$
as $u\rightarrow\infty$ where
$$
H({\bf T})=\Exp\exp(\sup_{t\in {\bf T}}\chi(t))<\infty\,.
$$
\end{lemma}

\begin{remark}
The continuity of the field $\chi(t)$ follows from Sudakov, Dudley and Fernique theorem (see \cite{pit:96}).
\end{remark}
\proof
\begin{eqnarray*}
\Prob(\sup_{t\in\, u^{-2/\alpha}{\bf T}}X(t)>u) &=&\frac{1}{\sqrt{2\pi}}\int_{-\infty}^{\infty}
e^{-\frac{v^2}{2}}\,\Prob(\sup_{t\in\, u^{-2/\alpha}{\bf T}}X(t)>u|X(0)=v)\,dv
\end{eqnarray*}
substituting $v=u-\frac{w}{u}$ 
$$
=\,\frac{1}{\sqrt{2\pi}u}\,e^{-\frac{u^2}{2}}\int_{-\infty}^\infty e^{w-\frac{w^2}{2u^2}}\,\Prob(\sup_{t\in\, u^{-2/\alpha}{\bf T}}X(t)>u|X(0)=u-\frac{w}{u})\,dw\,.
$$
Let us put 
$$
\chi_u(t)=u(X(u^{-2/\alpha}t)-u)+w\,.
$$
Thus let us rewrite the last integral without the function before the integral (which is $\Psi(u)$ as
$u\rightarrow\infty$)
\begin{equation}\label{meq1}
\int_{-\infty}^\infty e^{w-\frac{w^2}{2u^2}}\,\Prob(\sup_{t\in\, {\bf T}}\chi_u (t)>w|X(0)=u-\frac{w}{u})\,dw\,.
\end{equation}
Let us compute the expected value and variance of the distribution $\chi_u(t)$ under condition $X(0)=u-\frac{w}{u}$
(this distribution is Gaussian by Lemma \ref{gausv}). 
By Lemma \ref{gausv} we get
\begin{eqnarray*}
\Exp(\chi_u(t)|X(0))&=& u\Exp(X(u^{-2/\alpha}t)|X(0))-u^2+w\\
&=& u\alpha X(0)-u^2+w
\end{eqnarray*}
where $\alpha=r(u^{-2/\alpha}t)$. Hence
\begin{equation}\label{warexp}
ex(u,t)=\Exp(\chi_u(t)|X(0)=u-\frac{w}{u})=-u^2[1-r(u^{-2/\alpha}t)]+w[1-r(u^{-2/\alpha}t)]
\end{equation}
and by the assumptions it tends to $-|t_1|^\alpha-|t_2|^\alpha$ as $u\rightarrow\infty$.
Now let us calculate the variance
\begin{eqnarray}
\Var(\chi_u(t)|X(0)=u-\frac{w}{u})&=& u^2\Var(X(u^{-2/\alpha}t)|X(0)=u-\frac{w}{u})\nonumber\\
&=& u^2\Var(Z)\nonumber\\
&=& u^2(1-r^2(u^{-2/\alpha}t))\label{warvar}
\end{eqnarray}
where $Z$ in the second line is a suitable random variable from Lemma \ref{gausv} and 
by the assumptions it tends to $2(|t_1|^\alpha+|t_2|^\alpha)$ as $u\rightarrow\infty$.
Similarly we compute  
$$
\Var(\chi_u(t)-\chi_u(s)|X(0)=u-\frac{w}{u})= u^2\Var(X(u^{-2/\alpha}t)-X(u^{-2/\alpha}s)|X(0)=u-\frac{w}{u})
$$
by Lemma \ref{gausv}
$$
\,\,=\, u^2[\Var(X(u^{-2/\alpha}t)-X(u^{-2/\alpha}s))-[r(u^{-2/\alpha}t)-r(u^{-2/\alpha}s)]^2]\,.
$$
Thus we get
$$
\Var(\chi_u(t)-\chi_u(s)|X(0)=u-\frac{w}{u})=u^2[2[1-r(u^{-2/\alpha}(t-s))]-[r(u^{-2/\alpha}t)-r(u^{-2/\alpha}s)]^2]
$$
and one can estimate 
\begin{eqnarray*}
\Var(\chi_u(t)-\chi_u(s)|X(0)=u-\frac{w}{u})&\leq & 2u^2[1-r(u^{-2/\alpha}(t-s))]
\end{eqnarray*}
$$
\,\,\,\,\,= \,2(|t_1-s_1|^\alpha+|t_2-s_2|^\alpha)+u^2o(u^{-2}[|t_1-s_1|^\alpha+|t_2-s_2|^\alpha])
$$
$$
\,\,\,\,\,= \,(|t_1-s_1|^\alpha+|t_2-s_2|^\alpha)(2+o(1))
$$
where $o(1)\rightarrow 0$ if $u\rightarrow\infty$ or $|t_1-s_1|\rightarrow 0$ and $|t_2-s_2|\rightarrow 0$.
Hence
\begin{equation}\label{warcias}
\Var(\chi_u(t)-\chi_u(s)|X(0)=u-\frac{w}{u})\leq 3(|t_1-s_1|^\alpha+|t_2-s_2|^\alpha)
\end{equation}
for $u$ sufficiently large and $t,s$ belonging to a any bounded set of $\RL^2$.
One can also show that the covariance of $\chi_u(t)$ and $\chi_u(s)$ under condition $X(0)=u-\frac{w}{u}$
tends to $|t_1|^\alpha+|t_2|^\alpha+|s_1|^\alpha+|s_2|^\alpha-|t_1-s_1|^\alpha-|t_2-s_2|^\alpha$. Thus the finite dimensional distributions of the field $\chi_u(t)$ under condition $X(0)=u-\frac{w}{u}$ converge to
the finite dimensional distributions of $\chi(t)$ and by (\ref{warcias}) the distributions of the field $\chi_u(t)$ under condition $X(0)=u-\frac{w}{u}$ are tight which yield that the field $\chi_u(t)$ under condition $X(0)=u-\frac{w}{u}$ converges
weakly to $\chi(t)$ as $u\rightarrow\infty$.

From the weak convergence
\begin{equation}\label{zbsup}
\Prob(\sup_{t\in{\bf T}}\chi_u(t)>w|X(0)=u-\frac{w}{u})\rightarrow
\Prob(\sup_{t\in{\bf T}}\chi(t)>w)
\end{equation}
as $u\rightarrow\infty$.
Since the process $\chi_u(t)$ under condition $X(0)=u-\frac{w}{u}$ is continuous on ${\bf T}$
we get by Borell Theorem \ref{borel} that
$$
\Exp(\sup_{t\in{\bf T}}(\chi_u(t)-ex(u,t))|X(0)=u-\frac{w}{u})= m<\infty\,,
$$
$$
\sup_{t\in{\bf T}}\Var(\chi_u(t)|X(0)=u-\frac{w}{u})= \sigma^2<\infty
$$
where by (\ref{warexp}), (\ref{warvar}) and (\ref{zbsup}) $m$  and $\sigma^2$ depend only on $\alpha$
and
\begin{equation}\label{borelgl}
\Prob(\sup_{t\in{\bf T}}(\chi_u(t)-ex(u,t))>w|X(0)=u-\frac{w}{u})
\leq \exp\left(\frac{-(w-m)^2}{2\sigma^2}\right)
\end{equation}
for all $w\geq m$ for sufficiently large $u$. Since
$$
\Prob(\sup_{t\in{\bf T}}(\chi_u(t)-m)>w|X(0)=u-\frac{w}{u})\leq
\Prob(\sup_{t\in{\bf T}}(\chi_u(t)-ex(u,t))>w|X(0)=u-\frac{w}{u})
$$
and by (\ref{borelgl}) we have
\begin{equation}\label{borelgll}
\Prob(\sup_{t\in{\bf T}}\chi_u(t)>w|X(0)=u-\frac{w}{u})
\leq \exp\left(\frac{-(w-2m)^2}{2\sigma^2}\right)\,.
\end{equation}
Then using (\ref{borelgll}) the dominated convergence theorem yields
that
$$
\Exp [\exp(\sup_{t\in{\bf T}}\chi_u(t))|X(0)=u-\frac{w}{u}]\rightarrow
\Exp [\exp(\sup_{t\in{\bf T}}\chi(t))]
$$
as $u\rightarrow\infty$ and $\Exp [\exp(\sup_{t\in{\bf T}}\chi(t))]<\infty$. Thus taking into account 
(\ref{meq1}) we get the thesis.
\halmos

\begin{corollary}\label{estht}
If ${\bf T}=[a,b]\times [c,d]$ then
$$
H({\bf T})\leq \left\lceil b-a \right\rceil\left\lceil d-c \right\rceil
H([0,1]\times [0,1])
$$
where $\left\lceil x\right\rceil$ is the smallest integer larger than or equal to $x$. 
\end{corollary}
\proof 
We augment our rectangle to the rectangle with the sides of the length $\left\lceil b-a \right\rceil$ and
$\left\lceil d-c \right\rceil$. This rectangle can be divided into $\left\lceil b-a \right\rceil\left\lceil d-c \right\rceil$ unit squares. By the homogeneity of the random field $X$ we get the assertion. 
\halmos

Reducing one dimension in the previous lemma we get the following lemma.
\begin{lemma}\label{lemobok1}
Let $\chi(t)$ be a continuous stochastic Gaussian process where $t\in\RL$ with 
$\Exp\chi(t)=-|t|^\alpha$ and $\Cov(\chi(t),\chi(s))=|t|^\alpha+
|s|^\alpha-|t-s|^\alpha$ ($s\in\RL$) and $X(t)$ be a continuous stationary Gaussian 
process where $t\in\RL$ with expected value $\Exp X(t)=0$ and covariance 
$$
r(t)=\Exp(X(t+s)X(s))=1-|t|^\alpha+o(|t|^\alpha)\,.
$$
Then for any $T>0$
$$
\Prob(\sup_{t\in\, [0,\,u^{-2/\alpha}T]}X(t)>u)=\Psi(u)H(T)(1+o(1))
$$
as $u\rightarrow\infty$ where
\begin{equation}\label{defht}
H(T)=\Exp\exp(\sup_{t\in [0,\,T]}\chi(t))<\infty\,.
\end{equation}
\end{lemma}
\begin{remark}
Let us notice that $\chi(t)=B_H(t)-|t|^\alpha$ where $B_H$ is the fractional Brownian motion
with Hurst parameter $H=\alpha/2$ and $\Exp B_H^2(1)=2$. 
\end{remark}
\proof
The proof goes the same way as the proof of Lemma \ref{lemobok}.
\halmos

\begin{corollary}\label{estht1w}
For $T>0$
$$
H(T)\leq \left\lceil T \right\rceil
H([0,1])\,.
$$
\end{corollary}

The next lemma is different than Lemma D.2. in Piterbarg \cite{pit:96} that is the constant before 
exponent depends on $T$.
\begin{lemma}\label{lemdlak}
Let $0<\epsilon<1/2$ and $0<\epsilon^\alpha<1/2$ and $1-2|t|^\alpha\leq r(t)\leq 1-\frac{1}{2}|t|^\alpha$
for all $t\in [0,\epsilon]$ where $X(t)$ is defined in Lemma \ref{lemobok1}.
Then for $T>0$, $t_0>T$ and $u$ sufficiently large
$$
\Prob(\sup_{t\in\, [0,\,u^{-2/\alpha}T]}X(t)>u,\,\sup_{t\in\, [u^{-2/\alpha}t_0,\,u^{-2/\alpha}(t_0+T)]}X(t)>u)
\leq C(\alpha,t_0,T)\,\Psi(u)
$$
where
$$
C(\alpha,t_0,T)=
4\lceil D\,T \rceil
\,\lceil D\,(t_0+T) \rceil
\exp(-\frac{1}{8}(t_0-T)^\alpha)H([0,1]\times [0,1])\,.
$$
and  $D=\left(\frac{2\sqrt{2}}{\sqrt{7}}\right)^{2/\alpha}16^{1/\alpha}$.
\end{lemma}
\begin{remark}
Let us notice that the assumption $r(t)=1-|t|^\alpha+o(|t|^\alpha)$ implies that
there exists $\epsilon>0$ such that $1-2|t|^\alpha\leq r(t)\leq 1-\frac{1}{2}|t|^\alpha$ for all $t\in [0,\epsilon]$.
\end{remark}
\proof
Let us consider a Gaussian field $Y(t,s)=X(t)+X(s)$. Then
\begin{equation}\label{esty}
\Prob(\sup_{t\in A} X(t)>u,\,\sup_{t\in B} X(t)>u)
\leq
\Prob(\sup_{(t,s)\in\,A\times B} Y(t,s)>2u)
\end{equation}
where $A=[0,\,u^{-2/\alpha}T]$ and $B=[u^{-2/\alpha}t_0,\,u^{-2/\alpha}(t_0+T)]$.
Let us notice
\begin{eqnarray}
\sigma^2(t,s)&=& \Var Y(t,s)\label{vary}\\
&=& 2+2r(t-s)\nonumber\\
&=& 4-2(1-r(t-s))\nonumber\,.
\end{eqnarray}
From the assumptions of the lemma for $|t-s|\leq \epsilon$ we have
$$
\frac{1}{2}|t-s|^\alpha\leq 1-r(t-s)\leq 2|t-s|^\alpha 
$$
which gives
$$
4-4|t-s|^\alpha\leq \sigma^2(t,s)\leq 4-|t-s|^\alpha \,.
$$
Thus for sufficiently large $u$ we get 
\begin{equation}\label{infr}
\inf_{(t,s)\in (A\times B)}\sigma^2(t,s)\geq 4-4\sup_{(t,s)\in (A\times B)}|t-s|^\alpha
\geq 4-4\epsilon^\alpha> 2
\end{equation}
where in the last inequality we used the assumption of the lemma. Similarly
for sufficiently large $u$ we obtain
\begin{eqnarray}
\sup_{(t,s)\in (A\times B)}\sigma^2(t,s)&\leq&4-\inf_{(t,s)\in (A\times B)}|t-s|^\alpha\nonumber\\
&\leq& 4-|u^{-2/\alpha}(t_0-T)|^\alpha\nonumber\\
&=& 4-u^{-2}(t_0-T)^\alpha\,.\label{supsig}
\end{eqnarray}
Let us put 
$$
Y^{*}(t,s)=\frac{Y(t,s)}{\sigma(t,s)}
$$
where $\sigma(t,s)$ is defined in $(\ref{vary})$.
Let us estimate the right hand side of (\ref{esty}). Thus for sufficiently large $u$ we have
\begin{eqnarray}
\Prob(\sup_{(t,s)\in\,A\times B} Y(t,s)>2u)&=&
\Prob(\exists (t,s)\in\,A\times B: \frac{Y(t,s)}{\sigma(t,s)}>\frac{2u}{\sigma(t,s)})\nonumber\\
&\leq & \Prob(\sup_{(t,s)\in\,A\times B} Y^{*}(t,s)>\frac{2u}{\sqrt{4-u^{-2}(t_0-T)^\alpha}})
\label{supy}
\end{eqnarray}
where in the last line we used (\ref{supsig}).
Let us compute the following expectation for $(t,s)\in A\times B$ and $(t_1,s_1)\in A\times B$
\begin{eqnarray*}
\Exp[Y^{*}(t,s)-Y^{*}(t_1,s_1)]^2&=&\Exp\left[\frac{Y(t,s)-Y(t_1,s_1)}{\sigma(t,s)}+\frac{Y(t_1,s_1)}{\sigma(t,s)}
-\frac{Y(t_1,s_1)}{\sigma(t_1,s_1)}\right]^2\\
&\leq & 2\Exp\left[\frac{Y(t,s)-Y(t_1,s_1)}{\sigma(t,s)}\right]^2+\\
&&
\,\,\,\,\,\,\,2\left[\frac{1}{\sigma(t,s)}-\frac{1}{\sigma(t_1,s_1)}\right]^2
\Exp Y^2(t_1,s_1)
\end{eqnarray*}
where in the last inequality we used that $(a+b)^2\leq 2a^2+2b^2$ and continuing 
\begin{eqnarray*}
&\leq &\frac{2}{\inf_{(t,s)\in A\times B}\sigma^2(t,s)}\,\Exp\left[Y(t,s)-Y(t_1,s_1)\right]^2+\\
&&
\,\,\,\,\,\,\,\,\,\,\,\,\,\,\,\,\,\,2\left[\frac{1}{\sigma(t,s)}-\frac{1}{\sigma(t_1,s_1)}\right]^2
\sigma^2(t_1,s_1)\\
&=&
\frac{2}{\inf_{(t,s)\in A\times B}\sigma^2(t,s)}\,\Exp\left[Y(t,s)-Y(t_1,s_1)\right]^2+
2\left[\frac{\sigma(t_1,s_1)-\sigma(t,s)}{\sigma(t,s)}\right]^2\\
&\leq & \frac{2}{\inf_{(t,s)\in A\times B}\sigma^2(t,s)}\,
\left[\Exp\left[Y(t,s)-Y(t_1,s_1)\right]^2+[\sigma(t_1,s_1)-\sigma(t,s)]^2\right]
\end{eqnarray*}
using (\ref{infr}) for sufficiently large $u$ we get
\begin{eqnarray*}
&\leq & \Exp\left[Y(t,s)-Y(t_1,s_1)\right]^2+[\sigma(t_1,s_1)-\sigma(t,s)]^2\\
&=& \Exp[X(t)-X(t_1)+X(s)-X(s_1)]^2+[\sigma(t_1,s_1)-\sigma(t,s)]^2\\
&\leq & 2\Exp[X(t)-X(t_1)]^2+2\Exp[X(s)-X(s_1)]^2+[\sigma(t_1,s_1)-\sigma(t,s)]^2
\end{eqnarray*}
where in the last inequality we used that $(a+b)^2\leq 2a^2+2b^2$ and continuing 
\begin{eqnarray*}
&=& 2\Exp[X(t)-X(t_1)]^2+2\Exp[X(s)-X(s_1)]^2+\\
&&\,\,\,\,\,\,\,\,\,\,\,\,\,\,\,\sigma^2(t_1,s_1)-2\sigma(t_1,s_1)\sigma(t,s)+
\sigma^2(t,s)\\
&=& 2\Exp[X(t)-X(t_1)]^2+2\Exp[X(s)-X(s_1)]^2+\\
&&\,\,\,\,\,\,\,\,\,\,\,\,\,\,\,\Exp Y^2(t_1,s_1)-2\sqrt{\Exp Y^2(t_1,s_1)\Exp Y^2(t,s)}+
\Exp Y^2(t,s)
\end{eqnarray*}
by Schwarz inequality we obtain
\begin{eqnarray*}
&\leq & 2\Exp[X(t)-X(t_1)]^2+2\Exp[X(s)-X(s_1)]^2+\\
&&\,\,\,\,\,\,\,\,\,\,\,\,\,\,\,\Exp Y^2(t_1,s_1)-2\Exp[Y(t_1,s_1)Y(t,s)]+
\Exp Y^2(t,s)\\
&=& 2\Exp[X(t)-X(t_1)]^2+2\Exp[X(s)-X(s_1)]^2+\\
&&\,\,\,\,\,\,\,\,\,\,\,\,\,\,\,\Exp [Y(t,s)-Y(t_1,s_1)]^2\\
&=& 2\Exp[X(t)-X(t_1)]^2+2\Exp[X(s)-X(s_1)]^2+\\
&&\,\,\,\,\,\,\,\,\,\,\,\,\,\,\,\Exp[X(t)-X(t_1)+X(s)-X(s_1)]^2
\end{eqnarray*}
using the inequality $(a+b)^2\leq 2a^2+2b^2$ we get
\begin{eqnarray}
&\leq & 4\Exp[X(t)-X(t_1)]^2+4\Exp[X(s)-X(s_1)]^2\label{estexp}\,.
\end{eqnarray}
Since for $|t-t_1|\leq \epsilon$
\begin{eqnarray}
\Exp[X(t)-X(t_1)]^2&= & 2-2r(|t-t_1|)\nonumber\\
&\leq & 4|t-t_1|^\alpha\label{estexp1}
\end{eqnarray}
where in the last inequality we used the assumption of the lemma. Thus by (\ref{estexp}) and
(\ref{estexp1}) we have for $(t,s)\in A\times B$ and $(t_1,s_1)\in A\times B$ and $u$ sufficiently
large
\begin{equation}\label{estexp2}
\Exp[Y^{*}(t,s)-Y^{*}(t_1,s_1)]^2\leq 
16[|t-t_1|^\alpha+|s-s_1|^\alpha]\,.
\end{equation} 
Since $\Exp [Y^{*}(t,s)]^2=1$ and by (\ref{estexp2})
\begin{equation}\label{covstar}
\Exp[Y^{*}(t,s)Y^{*}(t_1,s_1)]\geq 1-8|t-t_1|^\alpha-8|s-s_1|^\alpha\,.
\end{equation}
Let us define the following random field
\begin{equation}\label{fieldz}
Z(t,s)=\frac{1}{\sqrt{2}}(\eta_1(t)+\eta_2(s))
\end{equation}
where $\eta_1$ and $\eta_2$ are independent Gaussian stationary processes with
$\Exp\eta_1(t)=\Exp\eta_2(t)=0$ and $\Exp[\eta_i(t)\eta_i(s)]=\exp(-32|t-s|^\alpha)$ for
$i=1,2$. Hence
\begin{eqnarray}
\Exp[Z(t,s)Z(t_1,s_1)]&=&\frac{1}{2}(\Exp[\eta_1(t)\eta_1(t_1)+\Exp[\eta_2(s)\eta_2(s_1)])\nonumber\\
&=& \frac{1}{2}[\exp(-32|t-t_1|^\alpha)+\exp(-32|s-s_1|^\alpha)]\nonumber\\
&\leq& 1-8|t-t_1|^\alpha-8|s-s_1|^{\alpha}\label{covz}
\end{eqnarray} 
for sufficiently small $|t-t_1|$ and $|s-s_1|$ by the fact that $e^{-x}\leq 1-\frac{1}{2}x$ for sufficiently small and positive $x$.
Thus by (\ref{covstar}) and (\ref{covz}) it follows
\begin{equation}\label{ineqcov}
\Exp[Y^{*}(t,s)Y^{*}(t_1,s_1)]\geq \Exp[Z(t,s)Z(t_1,s_1)]
\end{equation}
for sufficiently small $|t-t_1|$ and $|s-s_1|$. Hence by Slepian inequality we have for large $u$
\begin{equation}\label{supineq}
\Prob(\sup_{(t,s)\in A\times B}Y^{*}(t,s)>u^{*})\leq \Prob(\sup_{(t,s)\in A\times B}Z(t,s)>u^{*})
\end{equation}
where 
$$
u^{*}=\frac{2u}{\sqrt{4-u^{-2}(t_0-T)^\alpha}}
$$ 
(see (\ref{supy})).
Let us put
$$
\eta(t,s)=Z\left(\frac{t}{16^{1/\alpha}},\frac{s}{16^{1/\alpha}}\right)
$$
then
\begin{equation}\label{zeteta}
\Prob(\sup_{(t,s)\in\,A\times B} Z(t,s)>u^{*})
= \Prob(\sup_{(t,s)\in\,A'\times B'} \eta(t,s)>u^{*})
\end{equation}
where $A'=[0,\,u^{-2/\alpha}T 16^{1/\alpha}]$ and $B'=[u^{-2/\alpha}t_0 16^{1/\alpha},\,u^{-2/\alpha}(t_0+T)16^{1/\alpha}]$.
Let us notice that $\eta(t,s)$ satisfies the assumptions of Lemma \ref{lemobok} (for field $X$).
For 
$$
u\geq u_0=\left[\frac{(t_0-T)}{\epsilon}\right]^{\alpha/2}
$$
we get
$$
\frac{u^{*}}{u}=\frac{2}{\sqrt{4-u^{-2}(t_0-T)^\alpha}}
\leq  \frac{2}{\sqrt{4-u_0^{-2}(t_0-T)^\alpha}}
= \frac{2}{\sqrt{4-\epsilon^\alpha}}
<  \frac{2\sqrt{2}}{\sqrt{7}}
$$
where in the last inequality we used the assumption of the lemma that $\epsilon^\alpha<\frac{1}{2}$.
Thus it follows that $A'\subset[0,\,(u^{*}\frac{\sqrt{7}}{2\sqrt{2}})^{-2/\alpha}T 16^{1/\alpha}]$ and
$B'\subset [0,\,(u^{*}\frac{\sqrt{7}}{2\sqrt{2}})^{-2/\alpha}(t_0+T)16^{1/\alpha}]$. Let us define
${\bf T}=[0,\,(\frac{\sqrt{7}}{2\sqrt{2}})^{-2/\alpha}T 16^{1/\alpha}]\times [0,\,(\frac{\sqrt{7}}{2\sqrt{2}})^{-2/\alpha}(t_0+T)16^{1/\alpha}]$. Hence
\begin{eqnarray}
\Prob(\sup_{(t,s)\in\,A'\times B'} \eta(t,s)>u^{*})&\leq&
\Prob(\sup_{(t,s)\in\,(u^{*})^{-2/\alpha} {\bf T}} \eta(t,s)>u^{*})\nonumber\\
&=& \Psi(u^{*})H({\bf T})(1+o(1))\label{esteta}
\end{eqnarray}
as $u\rightarrow\infty$ where in the last line we used Lemma \ref{lemobok}.
By the fact that $\frac{1}{1-x}\geq 1+x$ for $x<1$ we get for sufficiently large $u$
$$
(u^{*})^2=\frac{4u^2}{4-u^{-2}(t_0-T)^\alpha}\geq
u^2[1+\frac{1}{4}u^{-2}(t_0-T)^\alpha]=u^2+\frac{1}{4}(t_0-T)^\alpha\geq u^2\,.
$$
Thus using (\ref{glasym}) we deduce that for sufficiently large $u$
$$
\Psi(u^*)\leq 2\Psi(u)\exp(-\frac{1}{8}(t_0-T)^\alpha)\,.
$$
Hence by (\ref{esteta}) it follows for sufficiently large $u$
\begin{eqnarray}
\Prob(\sup_{(t,s)\in\,A'\times B'} \eta(t,s)>u^{*})&\leq&
2\Psi(u)\exp(-\frac{1}{8}(t_0-T)^\alpha)H({\bf T})(1+o(1))\nonumber\\
&\leq &4\Psi(u)\exp(-\frac{1}{8}(t_0-T)^\alpha)H({\bf T})\label{etah}\,.
\end{eqnarray}
From Corollary \ref{estht} we obtain that
\begin{equation}\label{hest}
H({\bf T})\leq
H([0,1]\times [0,1])\lceil (\frac{\sqrt{7}}{2\sqrt{2}})^{-2/\alpha}T 16^{1/\alpha} \rceil
\lceil (\frac{\sqrt{7}}{2\sqrt{2}})^{-2/\alpha}(t_0+T)16^{1/\alpha} \rceil\,.
\end{equation}
Thus collecting (\ref{esty}), (\ref{supy}), (\ref{supineq}), (\ref{zeteta}), (\ref{etah})
and (\ref{hest}) we get the assertion of the lemma.
\halmos

\section{Pickands theorem}

\begin{theorem}
{\em (Pickands)} Let $X(t)$ where $t\in [0, p]$ be a continuous stationary Gaussian 
process with expected value $\Exp X(t)=0$ and covariance 
$$
r(t)=\Exp(X(t+s)X(s))=1-|t|^\alpha+o(|t|^\alpha)\,.
$$
Furthermore we assume that $r(t)<1$ for all $t>0$. Then
$$
\Prob(\sup_{t\in [0,p]} X(t)>u)=H_\alpha\, p\, u^{2/\alpha}\,\Psi(u)(1+o(1))
$$
as $u\rightarrow\infty$
where 
$$
H_{\alpha}=\lim_{T\rightarrow\infty}\frac{H(T)}{T}
$$
is positive and finite (Pickands constant) where $H(T)$ is defined in (\ref{defht}).
\end{theorem}
{\bf Proof:}
Put
$$
\Delta_k=[ku^{-2/\alpha}T, (k+1)u^{-2/\alpha}T]
$$
where $k\in \bbN$ and $T\geq p$ and $N_p=\left\lfloor \frac{p}{u^{-2/\alpha}T}\right\rfloor$.
Thus
\begin{eqnarray*}
\Prob(\sup_{t\in [0,p]} X(t)>u)&\leq& \sum_{k=0}^{N_p}\Prob(\sup_{t\in \Delta_k} X(t)>u)\\
&=& (N_p+1)\Prob(\sup_{t\in \Delta_0} X(t)>u)
\end{eqnarray*}
where in the last equality we use stationarity of the process $X$.
Thus using Lemma \ref{lemobok1} we get
\begin{equation}\label{nierpod}
\limsup_{u\rightarrow\infty}\frac{\Prob(\sup_{t\in [0,p]} X(t)>u)}{u^{2/\alpha}\Psi(u)}
\leq \frac{p}{T}\,H(T)\,.
\end{equation}
Let us estimate our probability from below 
\begin{eqnarray}
\Prob(\sup_{t\in [0,p]} X(t)>u)&\geq & \Prob(\bigcup_{k=0}^{N_p-1}\{\sup_{t\in \Delta_k} X(t)>u\})\nonumber\\
&\geq & N_p\,\Prob(\sup_{t\in \Delta_0} X(t)>u)\label{bonsup}\\
&& \,\,\,-\sum_{0\leq i<j\leq N_p-1}\Prob(\sup_{t\in \Delta_i} X(t)>u,\,\sup_{t\in \Delta_j} X(t)>u)\nonumber
\end{eqnarray}
where in the last inequality we applied Lemma \ref{bonf}. Let us consider the last double sum
(that is why the method is called double sum method)
\begin{eqnarray*}
\Sigma_2&=&\sum_{0\leq i<j\leq N_p-1}\Prob(\sup_{t\in \Delta_i} X(t)>u,\,\sup_{t\in \Delta_j} X(t)>u)\\
&=& \sum_{k=1}^{N_p-1}(N_p-k)\Prob(\sup_{t\in \Delta_0} X(t)>u,\,\sup_{t\in \Delta_k} X(t)>u)\\
&\leq & N_p\,\Prob(\sup_{t\in \Delta_0} X(t)>u,\,\sup_{t\in \Delta_1} X(t)>u)\\
&&\,\,\,\, +N_p\sum_{k=2}^{N_{\epsilon/4}-1}\Prob(\sup_{t\in \Delta_0} X(t)>u,\,\sup_{t\in \Delta_k} X(t)>u)\\
&&\,\,\,\, +N_p\sum_{k=N_{\epsilon/4}}^{N_{p}-1}\Prob(\sup_{t\in \Delta_0} X(t)>u,\,\sup_{t\in \Delta_k} X(t)>u)\,.
\end{eqnarray*}
Let us denote the last three terms by $A_1$, $A_2$ and $A_3$, respectively. 
We will show that these therms are negligible after dividing them by $u^{2/\alpha}\Psi(u)$ and passing with $u\rightarrow\infty$ and $T\rightarrow\infty$. Moreover bounds on them justify that Pickands constant is well-defined.

First let us consider $A_3$ and take $u$ such that $u^{-2/\alpha}T\leq \epsilon/16$. Then it is easy to notice that the distance of the intervals $\Delta_0$ and $\Delta_k$ is at least $\epsilon/4$ in $A_3$. Hence in $A_3$ (for $k$ from $A_3$) for $(t,s)\in \Delta_0\times\Delta_k$ we have 
\begin{eqnarray}
\Var(X(t)+X(s))&=& 2+2r(t-s)\nonumber\\
&=& 4-2(1-r(t-s))\nonumber\\
&\leq & 4-2\inf_{s\geq \epsilon/4}(1-r(s))\nonumber\\
&=& 4-\delta<4\label{borel1}
\end{eqnarray}
where $\delta=2\inf_{s\geq \epsilon/4}(1-r(s))>0$ (using the assumptions on $r(t)$).
Let us notice that $X(t)+X(s)$ is a continuous Gaussian field on $[0,T]\times [0,T]$
which implies by Borell Theorem \ref{borel} that
\begin{equation}\label{borel3}
\Exp\sup_{(t,s)\in \Delta_0\times\Delta_k}(X(t)+X(s))\leq m
\end{equation}
and by (\ref{borel1}) and (\ref{borel3}) we get
\begin{eqnarray*}
\Prob(\sup_{t\in \Delta_0} X(t)>u,\,\sup_{t\in \Delta_k} X(t)>u)&\leq&
\Prob(\sup_{(t,s)\in \Delta_0\times\Delta_k} X(t)+X(s)>2u)\\
&\leq &\exp\left(-\frac{(2u-m)^2}{2(4-\delta)}\right)\\
&=& \exp\left(-\frac{(u-m/2)^2}{2(1-\delta/4)}\right)\\
&\leq & \exp\left(-\frac{1}{2}\left(\frac{u-m/2}{1-\delta/8}\right)^2\right)
\end{eqnarray*}
where in the last inequality we used the fact that $1-\delta/4\leq (1-\delta/8)^2$.
Hence
\begin{eqnarray}
\limsup_{u\rightarrow\infty}\frac{A_3}{N_p \Psi(u)}
&\leq & \limsup_{u\rightarrow\infty}\frac{N_p^2\,\exp\left(-\frac{1}{2}\left(\frac{u-m/2}{1-\delta/8}\right)^2\right)}{N_p \Psi(u)}\nonumber\\
&=& \lim_{u\rightarrow\infty}\left\lfloor \frac{p}{u^{-2/\alpha}T}\right\rfloor
\sqrt{2\pi}\,u\,\exp(-\frac{1}{2}\left(\frac{u-a/2}{1-\delta/8}\right)^2
+\frac{1}{2}u^2)\nonumber\\
&=& 0 \label{lima3}
\end{eqnarray}
where the second line follows from (\ref{glasym}) and the fact that $1-\delta/8<1$ (by the assumption
$r(t)<1$ for $t>0$).

Now let us consider $A_2$. For $k\geq 2$ we have from Lemma \ref{lemdlak} ($C_1$ and $C_2$ constants depending on $\alpha$)
\begin{eqnarray*}
\lefteqn{\Prob(\sup_{t\in \Delta_0} X(t)>u,\,\sup_{t\in \Delta_k} X(t)>u)}\\
&\leq& C_1 \left\lceil C_2 T\right\rceil \left\lceil C_2 (k+1)T\right\rceil \exp(-\frac{1}{8}(k-1)^\alpha T^\alpha)
\Psi(u)\,.
\end{eqnarray*}
Thus 
\begin{eqnarray*}
A_2&\leq & C_1 \left\lceil C_2 T\right\rceil\Psi(u)N_p\sum_{k=2}^{N_{\epsilon/4}-1}\left\lceil C_2 (k+1)T\right\rceil \exp(-\frac{1}{8}(k-1)^\alpha T^\alpha)
\end{eqnarray*}
and let us estimate $\sum_{k=2}^{N_{\epsilon/4}-1}\left\lceil C_2 (k+1)T\right\rceil \exp(-\frac{1}{8}(k-1)^\alpha T^\alpha)$. We have
\begin{eqnarray*}
\lefteqn{\sum_{k=2}^{N_{\epsilon/4}-1}\left\lceil C_2 (k+1)T\right\rceil \exp(-\frac{1}{8}(k-1)^\alpha T^\alpha)}\\
&\leq& \sum_{k=2}^{\infty}\left\lceil C_2 (k+1)T\right\rceil \exp(-\frac{1}{8}(k-1)^\alpha T^\alpha)\\
&\leq &\left\lceil C_2T\right\rceil\sum_{k=2}^{\infty}(k+1) \exp(-\frac{1}{8}(k-1)^\alpha T^\alpha)\\
&=&\left\lceil C_2T\right\rceil\sum_{k=1}^{\infty}(k+2) \exp(-\frac{1}{8}k^\alpha T^\alpha)\\
&\leq& 3\left\lceil C_2T\right\rceil\sum_{k=1}^{\infty}k\exp(-\frac{1}{8}k^\alpha T^\alpha)\\
&\leq& 3\left\lceil C_2T\right\rceil\exp(-\frac{1}{8} T^\alpha)+3\left\lceil C_2T\right\rceil\int_{1}^\infty s\exp(-\frac{1}{8}s^\alpha T^\alpha)\,ds
\end{eqnarray*}
where the last inequality is valid for $T^\alpha>8/\alpha$ (then the function under integral is decreasing for $s>1$) and
substituting $t=\frac{1}{8}s^\alpha T^\alpha$ we continue (from now on $C$ will be any positive constant depending on $\alpha$ and its values can change from line to line)
\begin{eqnarray*}
&\leq& C\left\lceil T\right\rceil\exp(-\frac{1}{8} T^\alpha)+\frac{C\left\lceil T\right\rceil}{T^2}\int_{T^\alpha/8}^\infty t^{2/\alpha-1}\exp(-t)\,dt
\end{eqnarray*}
using the following property of the incomplete gamma function
$$
\int_u^\infty s^w e^{-s}\,ds = u^w e^{-u}(1+O(1/u))
$$
for $u\rightarrow\infty$ where $w\in \RL$ and
keeping on estimating we get
$$
\leq C\left\lceil T\right\rceil\exp(-\frac{1}{8}T^\alpha)(1+O(T^{-\alpha}))
$$
for $T^\alpha>8/\alpha$.
Thus we get
$$
A_2\leq C\left\lceil T\right\rceil^2\Psi(u)N_p\exp(-\frac{1}{8}T^\alpha)(1+O(T^{-\alpha}))
$$
which yields
\begin{equation}\label{lima2}
\limsup_{u\rightarrow\infty}\frac{A_2}{\Psi(u)N_p}\leq C\left\lceil T\right\rceil^2 \exp(-\frac{1}{8}T^\alpha)(1+O(T^{-\alpha}))\,.
\end{equation}
Now let us consider term $A_1$. Thus
\begin{eqnarray}
\lefteqn{\Prob(\sup_{t\in \Delta_0} X(t)>u,\,\sup_{t\in \Delta_1} X(t)>u)}\nonumber\\
&\leq&
\Prob(\sup_{t\in \Delta_0} X(t)>u,\,\sup_{t\in u^{-2/\alpha}[T,\,T+\sqrt{T}]} X(t)>u)\nonumber\\
&&+\Prob(\sup_{t\in \Delta_0} X(t)>u,\,\sup_{t\in u^{-2/\alpha}[T+\sqrt{T},\,2T+\sqrt{T}]} X(t)>u)\nonumber\\
&\leq& \Prob(\sup_{t\in u^{-2/\alpha}[T,\,T+\sqrt{T}]} X(t)>u)\nonumber\\
&&+\Prob(\sup_{t\in \Delta_0} X(t)>u,\,\sup_{t\in u^{-2/\alpha}[T+\sqrt{T},\,2T+\sqrt{T}]} X(t)>u)\nonumber\\
&=&\Prob(\sup_{t\in [0,\,u^{-2/\alpha}\sqrt{T}]} X(t)>u)\nonumber\\
&&+\Prob(\sup_{t\in \Delta_0} X(t)>u,\,\sup_{t\in u^{-2/\alpha}[T+\sqrt{T},\,2T+\sqrt{T}]} X(t)>u)\,.\label{2termA1}
\end{eqnarray}
First let us consider the second term of (\ref{2termA1}). By Lemma \ref{lemdlak} we have
$$
\Prob(\sup_{t\in \Delta_0} X(t)>u,\,\sup_{t\in u^{-2/\alpha}[T+\sqrt{T},\,2T+\sqrt{T}]} X(t)>u)
\,\,\,\,\,\,\,\,\,\,\,\,\,\,\,\,\,\,\,\,\,
$$
$$
\,\,\,\,\,\,\,\,\,\,\,\,\,\,\,\,
\leq 4\lceil C\,T \rceil
\,\lceil C\,(2T+\sqrt{T}) \rceil
\exp(-\frac{1}{8}\,T^{\alpha/2})H([0,1]\times [0,1])\Psi(u)\,.
$$
The first term from (\ref{2termA1}) can be estimated by Lemma (\ref{lemobok1})
$$
\Prob(\sup_{t\in [0,\,u^{-2/\alpha}\sqrt{T}]} X(t)>u)=
\Psi(u)H(\sqrt{T})(1+o(1))\,.
$$
Hence we obtain
\begin{eqnarray}
\lefteqn{\Prob(\sup_{t\in \Delta_0} X(t)>u,\,\sup_{t\in \Delta_1} X(t)>u)}\nonumber\\
&\leq&
\Psi(u)H(\sqrt{T})(1+o(1))\nonumber\\
&&+C\lceil T \rceil
\,\lceil 2T+\sqrt{T} \rceil
\exp(-\frac{1}{8}\,T^{\alpha/2})\Psi(u)\nonumber\\
&\leq& 
\Psi(u)\lceil \sqrt{T}\rceil H(1)(1+o(1))\nonumber\\
&&+C\lceil T \rceil
\,\lceil 2T+\sqrt{T} \rceil
\exp(-\frac{1}{8}\,T^{\alpha/2})\Psi(u)\label{ogra1}
\end{eqnarray}
where in the last inequality we used Corollary \ref{estht1w}. Thus we get
\begin{equation}\label{lima1}
\limsup_{u\rightarrow\infty}\frac{A_1}{N_p\Psi(u)}\leq
\lceil \sqrt{T}\rceil H(1)+C\lceil T \rceil
\,\lceil 2T+\sqrt{T} \rceil \exp(-\frac{1}{8}\,T^{\alpha/2})\,.
\end{equation}
Thus consider the lower bound
\begin{eqnarray*}
\liminf_{u\rightarrow\infty}\frac{\Prob(\sup_{t\in [0,p]} X(t)>u)}{p\,u^{2/\alpha}\Psi(u)}
&=& \liminf_{u\rightarrow\infty}\frac{\Prob(\sup_{t\in [0,p]} X(t)>u)}{N_p T\Psi(u)}
\end{eqnarray*}
which by Lemma \ref{lemobok1}, (\ref{bonsup}), (\ref{lima3}), (\ref{lima2}) and (\ref{lima1}) is bigger than or equal to
\begin{equation}\label{liminf}
f(T)=\frac{H(T)}{T}-
\frac{C\left\lceil T\right\rceil^2}{T}\exp(-\frac{1}{8}T^\alpha)(1+O(T^{-\alpha}))
\,\,\,\,\,\,\,\,\,\,\,\,\,\,\,\,\,\,\,\,\,\,\,\,
\end{equation}
$$
\,\,\,\,\,\,\,\,\,\,\,\,\,\,\,\,\,\,\,\,\,\,\,\,
-\frac{\lceil \sqrt{T}\rceil}{T} H(1)-C\frac{\lceil T \rceil}{T}
\,\lceil 2T+\sqrt{T} \rceil \exp(-\frac{1}{8}\,T^{\alpha/2})\,.
$$
Let us assume that $\limsup_{T\rightarrow\infty}\frac{H(T)}{T}>0$ then by (\ref{nierpod}) and
(\ref{liminf}) we get
\begin{eqnarray*}
\frac{H(T)}{T}&\geq &\limsup_{u\rightarrow\infty}\frac{\Prob(\sup_{t\in [0,1]} X(t)>u)}{u^{2/\alpha}\Psi(u)}\\
&\geq &\liminf_{u\rightarrow\infty}\frac{\Prob(\sup_{t\in [0,1]} X(t)>u)}{u^{2/\alpha}\Psi(u)}\\
&\geq & \limsup_{S\rightarrow\infty}f(S)\\
&=& \limsup_{S\rightarrow\infty}\frac{H(S)}{S}
\end{eqnarray*}
which implies
$$
\infty >\liminf_{T\rightarrow\infty}\frac{H(T)}{T}
\geq \limsup_{T\rightarrow\infty}\frac{H(T)}{T}>0
$$
and
$$
\lim_{T\rightarrow\infty}\frac{H(T)}{T}
$$
exists and is finite and positive.
It remains to prove that $\limsup_{T\rightarrow\infty}\frac{H(T)}{T}>0$.
Let us put $D=\bigcup_{j=0}^\infty\Delta_{2j}\cap [0,1]$. Then
$$
\Prob(\sup_{t\in [0,1]}X(t)>u)\geq \Prob(\sup_{t\in D}X(t)>u)\,.
$$
Applying Bonferroni inequality for the set $D$ (Lemma \ref{bonf} and see (\ref{bonsup}) 
and using Lemma \ref{lemobok1} and bound for $A_2$ and (\ref{lima3}) (note that $A_1$ disappears by the definition
of the set $D$) we get
\begin{eqnarray*}
\frac{H(T)}{T}&\geq &\limsup_{u\rightarrow\infty}\frac{\Prob(\sup_{t\in [0,1]} X(t)>u)}{u^{2/\alpha}\Psi(u)}\\
&\geq & \frac{H(S)}{2S}-\frac{C\left\lceil S\right\rceil^2}{S}\exp(-\frac{1}{8}S^\alpha)(1+O(S^{-\alpha}))\\
&= & S^{-1}(\frac{H(S)}{2}-C\left\lceil S\right\rceil^2\exp(-\frac{1}{8}S^\alpha)(1+O(S^{-\alpha})))
\end{eqnarray*}
which is positive for sufficiently large $S$ because $H(S)$ is increasing function of $S$ and
$C\left\lceil S\right\rceil^2\exp(-\frac{1}{8}S^\alpha)(1+O(S^{-\alpha}))$ tends to 0 when $S\rightarrow\infty$.
\halmos




\end{document}